\documentclass[uplatex,11pt,a4paper]{amsart} 

\usepackage{mathtools,amssymb,amsthm}
\usepackage{bm}
\usepackage{mathrsfs} 
\usepackage{stmaryrd} 

\usepackage{xcolor, graphicx, here}
\definecolor{darkgreen}{rgb}{0.0, 0.6, 0.0}

\usepackage{tikz}
\usetikzlibrary{intersections, calc, arrows, cd}

\usepackage[colorlinks,citecolor=darkgreen,linkcolor=red]{hyperref} 

\usepackage{geometry}
\numberwithin{equation}{section}
\numberwithin{figure}{section}
\numberwithin{table}{section}
\allowdisplaybreaks 

\usepackage{enumitem}
\setenumerate[1]{label={\upshape(\arabic*)},nosep}
\setenumerate[2]{label={\upshape(\alph*)},nosep}
\setitemize{nosep}
\setdescription{nosep}

\newlist{enumarabic}{enumerate}{1}
\setlist*[enumarabic]{label={\upshape(\arabic*)}, nosep}

\newlist{enumAlph}{enumerate}{1}
\setlist*[enumAlph]{label={\upshape(\Alph*)}, nosep}

\newlist{enumalph}{enumerate}{1}
\setlist*[enumalph]{label={\upshape(\alph*)}, nosep}

\newlist{enumRoman}{enumerate}{1}
\setlist*[enumRoman]{label={\upshape(\Roman*)}), nosep}

\newlist{enumroman}{enumerate}{1}
\setlist*[enumroman]{label={\upshape(\roman*)}, nosep}

\usepackage{cleveref}
 \crefname{section}{\S\!}{\S\S\!}
 \crefname{subsection}{\S\!}{\S\S\!}
 \crefname{subsubsection}{\S\!}{\S\S\!}

\theoremstyle{plain}
\newtheorem{theorem}{Theorem}[section]
 \crefname{theorem}{Theorem}{Theorems}
\newtheorem{theorem-definition}[theorem]{Theorem and Definition}
 \crefname{theorem-definition}{Theorem and Definition}{Theorem and Definition}
\newtheorem{lemma}[theorem]{Lemma}
 \crefname{lemma}{Lemma}{Lemmas}
\newtheorem{proposition}[theorem]{Proposition}
 \crefname{proposition}{Proposition}{Propositions}
\newtheorem{proposition-definition}[theorem]{Proposition and Definition}
 \crefname{proposition-definition}{Proposition and Definition}{Proposition and Definition}
\newtheorem{corollary}[theorem]{Corollary}
 \crefname{corollary}{Corollary}{Corollaries}

 \crefname{fact}{Fact}{Facts}
\newtheorem{claim}[theorem]{Claim}
 \crefname{claim}{Claim}{Claims}

 \crefname{conjecture}{Conjecture}{Conjectures}

 \crefname{problem}{Problem}{Problems}

 \crefname{question}{Question}{Questions}
\newtheorem{theoremA}{Theorem}
 \crefname{theoremA}{Theorem}{Theorems} 

\theoremstyle{definition}
\newtheorem{definition}[theorem]{Definition}
 \crefname{definition}{Definition}{Definitions}

 \crefname{notation}{Notation}{Notations}
\newtheorem{remark}[theorem]{Remark}
 \crefname{remark}{Remark}{Remarks}
\newtheorem{remark-definition}[theorem]{Remark and Definition}
 \crefname{remark-definition}{Remark and Definition}{Remark and Definition}
\newtheorem{example}[theorem]{Example}
 \crefname{example}{Example}{Examples}

 \crefname{observation}{Observation}{Observations}

 \crefname{construction}{Construction}{Constructions}

 \crefname{condition}{Condition}{Conditions}

\newenvironment{pfclaim}[1][$\Diamond$]
{\def\claimQED{{#1}}\noindent {\em Proof of the claim. }}
{\leavevmode\unskip\penalty9999 \hbox{}\nobreak\hfill
    \quad\hbox{\claimQED}{\smallskip}}

\newcommand{\ol}{\overline}
\newcommand{\ul}{\underline}

\newcommand{\wt}{\widetilde}

\newcommand{\xr}{\xrightarrow}

\newcommand{\C}{\mathbb{C}}
\newcommand{\F}{\mathbb{F}}

\newcommand{\Z}{\mathbb{Z}}

\newcommand{\G}{\mathbb{G}}

\newcommand{\m}{\mathfrak{m}}
\newcommand{\n}{\mathfrak{n}}

\newcommand{\scrO}{\mathscr{O}}

\DeclareMathOperator{\GL}{GL}
\DeclareMathOperator{\SL}{SL}

\DeclareMathOperator{\Spec}{Spec}

\newcommand{\id}{\mathrm{id}}

\DeclareMathOperator{\Hom}{Hom}

\newenvironment{psmatrix}{\left(\begin{smallmatrix}}{\end{smallmatrix}\right)}

\newcommand{\relmiddle}[1]{\mathrel{}\middle#1\mathrel{}} 

\begin{document}
\title[Frobenius lifts on invariant subrings]{Frobenius lifts on invariant subrings, with applications to rational double points and perfectoid towers}
\author[K.~Hayashi]{Kazuki Hayashi}
\address{Department of Mathematics, Institute of Science Tokyo, 2-12-1 Ookayama, Meguro, Tokyo 152-8551}
\email{hazuki0694@gmail.com}


\begin{abstract}
We study Frobenius lifts on invariant subrings and give applications to rational double points and perfectoid towers.
We realize rational double points of ADE over $W(k)$ as invariant subrings $A^G$, where $A$ is a ramified two-dimensional complete regular local ring and $G$ is a finite linearly reductive subgroup scheme of $\SL_{2,W(k)}$.
We also prove a general theorem on the existence of Frobenius lifts on invariant subrings: if a flat affine group scheme $G$ over $W(k)$ admits a functorial Reynolds operator and a relative Frobenius lift, then any formal power series ring equipped with a linear action of $G$ admits a Frobenius lift that restricts to a Frobenius lift on the invariant subring.
Building on these results, we construct perfectoid towers arising from rational double points in mixed characteristic. We identify their tilts with the perfect towers arising from rational double points in positive characteristic of the same ADE type.
\end{abstract}

\maketitle

\setcounter{tocdepth}{2}
\tableofcontents
\section{Introduction}

\emph{Rational double points} form one of the most fundamental classes of surface singularities. Over the complex numbers $\C$, they were first studied from several different points of view by F.~Klein, P.~du Val, and others. They are also known as du Val singularities, Kleinian singularities, or simple surface singularities. One of their most remarkable features is the appearance of the simply laced Dynkin diagrams $A_n$, $D_n$, $E_6$, $E_7$, and $E_8$. More precisely, the exceptional curves in the minimal resolution of a rational double point are smooth rational curves of self-intersection $-2$, and their dual intersection graph is one of the above Dynkin diagrams. Rational double points also arise as quotient singularities $\C^2/G$, where $G$ is a finite subgroup of $\SL_2(\C)$.  The correspondence between such finite subgroups and the ADE Dynkin diagrams is usually referred to as the Klein correspondence.

The modern algebraic theory of rational surface singularities was developed by M.~Artin \cite{Art66} and J.~Lipman \cite{Lip69}. Let $(R,\m,k)$ be a two-dimensional excellent normal local ring, and let $\pi\colon Y\to \Spec(R)$ be a resolution of singularities. The singularity is called \emph{rational} if $R^1\pi_*\scrO_Y=0$. A \emph{rational double point} is a two-dimensional Gorenstein rational singularity. If $R$ is complete, it is a hypersurface and can be written as $R\cong S/(f)$, where $(S,\n,k)$ is a three-dimensional complete regular local ring and $f\in\n^2\setminus\n^3$. In particular, $R$ has multiplicity two, which explains the terminology ``double point.''  Equivalently, in the classical geometric setting, rational double points are precisely the rational surface singularities whose exceptional curves in the minimal resolution all have self-intersection $-2$.

Over an algebraically closed field of characteristic zero, the theory is essentially the same as over $\C$. In particular, rational double points are classified by the ADE Dynkin diagrams and, after a suitable choice of formal coordinates, are defined by the familiar ADE equations. They can also be described as quotient singularities by finite subgroups of $\SL_2$. Thus, their defining equations, resolution graphs, and associated finite groups are related by the classical ADE correspondence.

The situation over an algebraically closed field $k$ of positive characteristic $p>0$ is more subtle. In sufficiently large characteristic, the characteristic-zero classification remains valid in essentially the same form. In small characteristics, however, additional phenomena occur: a fixed ADE type may contain more than one isomorphism class, the standard defining equations may require modifications, and the quotient description by a finite constant subgroup $G$ of $\SL_2(k)$ may fail, particularly when $p$ divides the order of $G$. Artin \cite{Art77} carried out a detailed study of rational double points in positive characteristic and proved, among other results, that every such singularity admits a finite, possibly ramified, cover by a regular scheme.
Nevertheless, an analogue of the Klein correspondence can be recovered in positive characteristic by replacing finite groups with finite \emph{group schemes}. Indeed, the two-dimensional $F$-rational Gorenstein complete local rings over an algebraically closed field of positive characteristic are classified by the ADE Dynkin diagrams, with the characteristic restrictions given in \cite[Example~18]{HL02} (see also \cite{WY00}). Building on this classification, M.~Hashimoto \cite{Has15} classified the linearly reductive finite subgroup schemes of $\SL_{2,k}$ up to conjugation and established a one-to-one correspondence between their conjugacy classes and the isomorphism classes of two-dimensional $F$-rational Gorenstein complete local rings with coefficient field $k$.
In particular, every rational double point in the $F$-rational range can be realized as a quotient singularity by a linearly reductive finite subgroup scheme of $\SL_{2,k}$. Thus Hashimoto's result provides an analogue of the classical Klein correspondence in positive characteristic.

It is natural to ask to what extent this picture extends to \emph{mixed characteristic}. For residual characteristic greater than $5$, J.~Carvajal-Rojas, L.~Ma, T.~Polstra, K.~Schwede, and K.~Tucker classified rational double points in mixed characteristic by showing that they admit the expected ADE normal forms. They also constructed a split finite cover by a regular scheme:

\begin{theorem}[{\cite{CRMP+21}}]
\label{thm:CRMP}
Let $(R,\m,k)$ be a complete local ring of a rational double point over $W(k)$, where $k$ is an algebraically closed field of characteristic $p\geq 7$.
  \begin{enumerate}
  \item Write $R\cong S/(f)$, where $(S,\n,k)$ is a three-dimensional complete regular local ring of mixed characteristic $(0,p)$ and $f\in\n^2\setminus\n^3$. Then there exists a choice of a regular system of parameters $x,y,z$ of $S$ such that, up to multiplication by a unit, $f$ can be written in one of the following forms:
\begin{equation}
\left\{ \,
    \begin{aligned}
    & x^2+y^2+z^{n+1} && A_n && (n\geq 1), \\
    & x^2+y^2z+z^{n-1} && D_n && (n\geq 4), \\
    & x^2+y^3+z^4 && E_6, && \\
    & x^2+y^3+yz^3 && E_7, && \\
    & x^2+y^3+z^5 && E_8. &&
    \end{aligned}
\right.
\end{equation}
  \item There exists a finite cover $\pi\colon Y\to X=\Spec(R)$ such that $\scrO_X\to \pi_*\scrO_Y$ splits as an $\scrO_X$-linear map and $Y$ is a regular scheme.
  \end{enumerate}
\end{theorem}

The theorem provides mixed-characteristic analogues of both the classical ADE classification and Artin's theorem on regular covers in positive characteristic. However, unlike the classical equal-characteristic description, it does not directly realize these singularities as invariant subrings of regular local rings. Our first main result supplies such a description (\cref{s:RDPG}). Namely, we show that every rational double point under consideration arises from a linear action of a linearly reductive finite subgroup scheme of $\SL_{2,W(k)}$ on a suitable ramified regular local ring.

\begin{theoremA}[\cref{thm:RDPG}]
\label{thm:A}
Let $(R,\m,k)$ be a complete local ring of a rational double point over $W(k)$, where $k$ is an algebraically closed field of characteristic $p>0$. Suppose that $R$ is of one of the ADE types listed in \textup{Table~\ref{tb:RDP}} and that $p$ satisfies the corresponding condition in the third column of the table.
Then there exist a finite linearly reductive subgroup scheme $G\subset\SL_{2,W(k)}$ of the type specified in the fourth column and a ramified two-dimensional complete regular local ring $A$ equipped with a $G$-action such that $R\cong A^G$.
\end{theoremA}

\begin{table}[htbp]
  \centering
  \caption{Rational double points as invariant subrings}
  \label{tb:RDP}
\begin{tabular}{ccccc} \hline
Type & Equation & Residual characteristic & $G$ & $\#G$ \\ \noalign{\hrule height 1.2pt}
$A_n\ (n\geq 1)$
  & $x^2+y^2+z^{n+1}$
  & $p\geq 3$ 
  & $\mu_{n+1}$
  & $n+1$
\\
$D_n\ (n\geq 4)$
  & $x^2+y^2z+z^{n-1}$
  & $p\geq 3$
  & $D_{4(n-2)}$
  & $4(n-2)$
\\
$E_6$
  & $x^2+y^3+z^4$
  & $p\geq 5$
  & $T$
  & $24$
\\
$E_7$
  & $x^2+y^3+yz^3$
  & $p\geq 5$
  & $O$
  & $48$
\\
$E_8$
  & $x^2+y^3+z^5$
  & $p\geq 7$
  & $I$
  & $120$ \\ \hline
\end{tabular}
\end{table}

The theorem complements the classification in \cref{thm:CRMP} (1) by giving a realization as invariant subrings. Moreover, it gives another proof of \cref{thm:CRMP} (2): since $G$ is finite linearly reductive, the inclusion $A^G\hookrightarrow A$ splits as an $A^G$-linear map (see \cref{lem:AGpure}).

We next study \emph{Frobenius lifts} on invariant subrings in mixed characteristic (\cref{s:FrobG}).
Let $W(k)$ be the Witt ring of a field $k$ of characteristic $p>0$. Then we have the Witt Frobenius $F\colon W(k)\to W(k)$, which is a Frobenius lift. The formal power series ring $A=W(k)\llbracket x_1,\ldots,x_n\rrbracket$ also admits a Frobenius lift $\varphi\colon A\to A$ defined by $\varphi|_{W(k)}=F$ and $\varphi(x_i)=x_i^p$ ($1\leq {}^\forall i\leq n$).
In general, however, this Frobenius lift need not preserve the invariant subring $A^G$. We therefore investigate when the standard Frobenius lift restricts to $A^G$, and, more generally, how to construct a Frobenius lift on $A$ that restricts to $A^G$. Our second main result gives a general criterion for such a construction.

\begin{theoremA}[{Special case of \cref{thm:GequivF}}]
\label{thm:B}
Let $G$ be a flat affine group scheme over the Witt ring $W(k)$ of a field $k$ of characteristic $p>0$.
Suppose that $G$ admits a functorial Reynolds operator \textup{(\cref{def:Reynolds})} and a relative Frobenius lift $F_{G/W(k)}\colon G\to G^{(p)}\coloneqq G\otimes_{W(k),F}W(k)$ \textup{(\cref{def:relFlift})}. Set $A\coloneqq W(k)\llbracket x_1,\ldots,x_n\rrbracket$ equipped with a linear action of $G$.
Then there exists a Frobenius lift $\varphi_{A,G}\colon A\to A$ restricts to a Frobenius lift on $A^G$.
\end{theoremA}

The existence of Frobenius lifts is closely related to the theory of \emph{perfectoid towers} introduced by S.~Ishiro, K.~Nakazato, and K.~Shimomoto \cite{INS25}. The \emph{tilting operation} $S\leadsto S^\flat$ for perfectoid rings $S$ makes bridge between mixed and positive characteristic. However, perfectoid rings arising in this context are typically non-Noetherian and therefore lie outside the usual framework of Noetherian singularity theory. A perfectoid tower
\[
R_0\to R_1\to \cdots \to R_i\to \cdots
\]
may be regarded as a Noetherian approximation of a perfectoid ring that retains a tilting operation at the level of a tower of Noetherian rings. It provides a perfect(oid) tower in positive characteristic
\[
R_0^{s.\flat}\to R_1^{s.\flat}\to \cdots \to R_i^{s.\flat}\to \cdots
\]
In particular, the \emph{$i$-th small tilt} $R_i^{s.\flat}$ makes it possible to relate a Noetherian ring of mixed
characteristic to one of positive characteristic. We refer to \cite{Ha26a} for further study of small tilts; see also \cite{HIS26, Ha26b,Ha26c}.
R.~Ishizuka \cite{Ish26} showed that a large class of Frobenius lifts gives rise to perfectoid towers through a construction from prisms; see also \cite{IS25,Ya26} for related works. We apply this construction to the invariant subrings concerning in \cref{thm:A} (\cref{s:RDPtower}). More precisely, by applying \cref{thm:B} to the group schemes and their actions in \cref{thm:A}, we construct Frobenius lifts on $W(k)\llbracket u,v\rrbracket^G$ and use them to obtain perfectoid towers.
Furthermore, the tilt of this tower arises from the invariant subring $k\llbracket u,v\rrbracket^{G_k}$ associated with the special fiber $G_k=G\otimes_{W(k)}k$ and is therefore a rational double point of the same ADE type as $R$:

\begin{theoremA}[\cref{thm:RDPtower}]
\label{thm:C}
Let $(R,\m,k)$ be the local ring of a rational double point over $W(k)$ as in \textup{\cref{thm:A}}. Then there exists a perfectoid tower arising from the pair $(R,(p))$
\[
R=R_0\to R_1\to \cdots \to R_i\to \cdots.
\]
Moreover, its tilt is isomorphic to
\[
R^{s.\flat}=k\llbracket u,v\rrbracket^{G_k} \xr{\varphi} k\llbracket u,v\rrbracket^{G_k} \xr{\varphi}\cdots \xr{\varphi} k\llbracket u,v\rrbracket^{G_k} \xr{\varphi}\cdots,
\]
where $\varphi$ is the absolute Frobenius. In particular, the $0$-th small tilt $R^{s.\flat}$ is a rational double point over $k$ of the same ADE type as $R$.
\end{theoremA}

Therefore, the perfectoid tower constructed in the theorem does more than provide a perfectoid approximation of $R$: its small tilt recovers the corresponding rational double point in positive characteristic. In this way, our construction gives a direct connection between mixed-characteristic rational double points and their positive-characteristic counterparts while preserving their ADE types. This compatibility suggests that perfectoid towers provide a useful framework for singularity theory between mixed and positive characteristic.

\medskip\noindent
\textbf{Acknowledgement.} 
The author would like to express sincere gratitude to Ryo Ishizuka, Ryu Tomonaga, Tatsuki Yamaguchi, and Shou Yoshikawa for helpful comments, discussions, and suggestions.

\medskip\noindent
\textbf{Convention.}
All rings are assumed to be commutative and unital.

\section{Preliminaries}

Throughout this section, $R$ is a ring and $G$ is a flat affine group scheme over $R$. We denote by $R[G]$ the coordinate ring of $G$, which is a commutative Hopf algebra over $R$.

\begin{definition}
  \begin{enumerate}
  \item A \emph{$G$-algebra} is an $R$-algebra $A$ together with an $R$-algebra coaction $\omega_A\colon A\to A\otimes_RR[G]$.
  \item Let $A$ be a $G$-algebra. A \emph{$(G,A)$-module} is an $A$-module $M$ together with a coaction $M\to M\otimes_RR[G]$ such that $\omega_M(ax)=\omega_A(a)\omega_M(x)$ for all $a\in A$ and $x\in M$. A $(G,R)$-module is simply called a \emph{$G$-module}.
  \item Let $A$ be a $G$-algebra, and $M$ a $(G,A)$-module. We define the \emph{invariant submodule} by
  \[
  M^G\coloneqq\{x\in M\mid \omega_M(x)=x\otimes 1\}.
  \]
  Note that $A^G$ is an $R$-subalgebra of $A$ and $M^G$ is an $A^G$-submodule of $M$.
  \item If $M,N$ are $G$-modules, then we say that an $R$-linear map $f\colon M\to N$ is \emph{$G$-equivariant} if the following diagram commutes:
  \[
  \begin{tikzcd}
  M \rar["\omega_M"] \dar["f"'] & {M\otimes_RR[G]} \dar["{f\otimes_RR[G]}"] \\
  N \rar["\omega_N"'] & {N\otimes_RR[G]}.
  \end{tikzcd}
  \]
  \end{enumerate}
\end{definition}

\begin{remark}
\label{rem:HomG}
Let $A$ be a $G$-algebra, and let $M$ and $N$ be $(G,A)$-modules such that $M$ is finitely presented as an $A$-module.
Then $\Hom_A(M,N)$ carries a canonical $(G,A)$-module structure. On valued points, the action is given by $(g\cdot f)(x)=g\cdot f(g^{-1}\cdot x)$ for $g\in G$, $f\in\Hom_R(M,N)$, and $x\in M$. We deduce from this the following.
  \begin{itemize}
  \item An element $f\in\Hom_R(M,N)$ is $G$-invariant if and only if the $R$-linear map $f\colon M\to N$ is $G$-equivariant.
  \end{itemize}
\end{remark}

The category of $G$-modules, whose morphisms are the $G$-equivariant $R$-linear maps, is a Grothendieck category due to the flatness of $G$ (\cite[Corollary 5.3.11]{Has}). We therefore make the following definition.

\begin{definition}[{\cite[Definition 3.1]{Wan25}}]
We say that $G$ is \emph{linearly reductive} if the invariants functor $(-)^G$ from the category of $G$-modules to the category of $R$-modules is exact.
\end{definition}

When $R$ is a field, there are several characterizations of linearly reductive group schemes.
One particularly useful characterization is in terms of Reynolds operators. We recall the following definitions and results from \cite{Wan25}, which apply over general base rings.

\begin{definition}[{\cite[Definition 3.3]{Wan25}}]
\label{def:Reynolds}
  \begin{enumerate}
  \item A \emph{Reynolds operator} on a $G$-module $M$ is an $R^G$-linear map $\rho_M\colon M\to M^G$ such that $\rho_M|_{M^G}=\id_{M^G}$.
  \item We say that \emph{$G$ admits a functorial Reynolds operator} if there exists a Reynolds operator $\rho_M$ for every $G$-module $M$, functorially in $M$.
  \end{enumerate}
\end{definition}

\begin{theorem}[{\cite[Propositions 3.4 and 3.6, Theorem 3.8]{Wan25}}]
\label{thm:Reynolds}
Consider the following conditions:
\begin{enumerate}
\item $G$ admits a functorial Reynolds operator.
\item There exists a Reynolds operator on $R[G]$.
\item $G$ is linearly reductive.
\end{enumerate}
Then the implications ``\textup{(1)} $\Leftrightarrow$ \textup{(2)} $\Rightarrow$ \textup{(3)}'' hold. If, moreover, $R$ is a Dedekind domain and $G$ is smooth over $R$, then all conditions are equivalent.
\end{theorem}

Whenever there exists a Reynolds operator on a $G$-module $M$, the canonical injection $M^G\hookrightarrow M$ splits as a map of $R^G$-modules for every $G$-module $M$. In particular, the injection is \emph{pure} as a map of $R^G$-modules, that is, the canonical map $M^G\otimes_{R^G}N\longrightarrow M\otimes_{R^G}N$ is injective for every $R^G$-module $N$.
The following lemma shows that the purity statement holds for any linearly reductive group scheme, without Reynolds operators.

\begin{lemma}
\label{lem:AGpure}
If $G$ is linearly reductive, then for any $G$-algebra $A$ and any $(G,A)$-module $M$, the inclusion $M^G\hookrightarrow M$ is pure as an $A^G$-linear map. In particular, for any ideal $I\subset A^G$, the inclusion $A^G\hookrightarrow A$ induces an isomorphism
\[
A^G/I\xr{\cong} (A/IA)^G.
\]
\end{lemma}

\begin{proof}
Let $N$ be an $A^G$-module. Then the injectivity of $M^G\otimes_{A^G}N\to M\otimes_{A^G}N$ follows from the canonical isomorphism
\[
M^G\otimes_{A^G}N \xr{\cong}(M\otimes_{A^G}N)^G.
\]
As remarked in \cite[Remark 4.8]{Alp13}, this isomorphism is a special case of the projection formula for the quotient stack $[\Spec(A)/G]\to\Spec(A^G)$ (\cite[Proposition 4.7]{Alp13}). In the present setting, this also follows by a simpler argument: as $N$ is an inductive limit of finitely presented $A^G$-modules, we may assume that $N$ is finitely presented. Taking a finite presentation of $N$ and using the exactness of the functor $(-)^G$, we reduce to the case where $N$ is a finite free $A^G$-module, in which case the assertion is immediate.
\end{proof}

Let us give some examples relevant to this paper.

\begin{example}
\label{ex:linred}
  \begin{enumerate}
  \item For a finite group $\Gamma$, let $G=\ul{\Gamma}_R$ be the associated constant group scheme. Then $G$ is linearly reductive if and only if $\#\Gamma\in R^\times$. In this case, a functorial Reynolds operator is given by $\rho_M(x)=\tfrac{1}{\#\Gamma}\sum\limits_{\gamma\in\Gamma}\gamma\cdot x$ for every $G$-module $M$.
  \item For an abelian group $\Delta$, let $G=D_R(\Delta)$ be the associated \emph{diagonalizable group scheme}, i.e., $D_R(\Delta)(X)$ is the set of group homomorphisms $\Delta\to\Gamma(X,\scrO_X)^\times$. The underlying $R$-scheme is $\Spec(R\Delta)$, where $R\Delta$ denotes the group algebra (\cite[Expos\'{e} I, 4.4]{SGA3-1}). Then a $G$-module is nothing but a $\Delta$-graded $R$-module $M=\bigoplus_{d\in\Delta}M_d$, and since $M^G=M_0$, one sees that $G$ is linearly reductive. Moreover, the canonical projection $\rho_M\colon M\to M_0=M^G$ defines a functorial Reynolds operator. In particular, the multiplicative group $\G_{m,R}=D_R(\Z)$ and the group of $n$-th roots of unity $\mu_{n,R}=D_R(\Z/n\Z)$ are linearly reductive and admit functorial Reynolds operators.
  \item Let $1\to N\to G\to H\to 1$ be an exact sequence of flat affine group schemes over $R$. Then $G$ is linearly reductive if and only if $N$ and $H$ are linearly reductive. One can prove it in a similar way to the proof of \cite[Lemma 2.2]{Has15}. Moreover, if both $N$ and $H$ admit functorial Reynolds operators, then so does $G$. Indeed, for any $G$-module $M$, $M^N$ is naturally an $H$-module, and we have $M^G=(M^N)^H$. Thus the composite $M\xr{\rho_M} M^N \xr{\rho_{M^N}} (M^N)^H=M^G$ gives a functorial Reynolds operator for $G$.
  \end{enumerate}
\end{example}

\section{Rational double points as invariant subrings}
\label{s:RDPG}

Throughout this section, $k$ is an algebraically closed field of characteristic $p>0$ and $W(k)$ is the ring of Witt vectors. All group schemes are over $W(k)$.

Let $(R,\m,k)$ be a two-dimensional complete Gorenstein rational singularity, and suppose that $R$ is not regular. Then $R$ is a rational double point, i.e., $R\cong S/(f)$, where $(S,\n,k)$ is a three-dimensional regular local ring and $f\in\n^2\setminus\n^3$ (\cite[Lemma 4.1]{CRMP+21}). Moreover, if $p\geq 7$, the possible defining equations $f$ are classified into the types $A_n$ ($n\geq 1$), $D_n$ ($n\geq 4$), $E_6$, $E_7$, and $E_8$; see \cite[Theorem B]{CRMP+21}.
Our main interest lies in hypersurface singularities of the above ADE types. Since the classification theorem requires the assumption $p\geq 7$, we define these types directly in terms of their defining equations, independently of the classification.

\begin{definition}
An \emph{ADE singularity} is a two-dimensional hypersurface $R$ for which there exist a three-dimensional complete regular local ring $(S,\n,k)$ of residual characteristic $p>0$, a regular system of parameters $x,y,z$ of $S$, and an isomorphism $R\cong S/(f)$, where $f$ and $p$ satisfy one of the conditions in the following table. 
\begin{table}[htbp]
  \centering
  \caption{ADE singularities}
  \label{tb:ADE}
\begin{tabular}{ccc}
\hline Type & Equation & Residual characteristic \\ \noalign{\hrule height 1.2pt}
$A_n$ ($n\geq 1$) & $x^2+y^2+z^{n+1}$ & $p\geq 3$ \\
$D_n$ ($n\geq 4$) & $x^2+y^2z+z^{n-1}$ & $p\geq 3$ \\
$E_6$ & $x^2+y^3+z^4$ & $p\geq 5$ \\
$E_7$ & $x^2+y^3+yz^3$ & $p\geq 5$ \\
$E_8$ & $x^2+y^3+z^5$ & $p\geq 7$
\\ \hline
\end{tabular}
\end{table}
\end{definition}

Under the indicated assumptions on the residual characteristic, a sequence of blow-ups yields the minimal resolution of each ADE singularity. Its dual graph is the corresponding ADE Dynkin diagram. It follows from Artin's rationality criterion that every ADE singularity is a rational singularity. Being Gorenstein, it is a rational double point. We therefore also refer to an ADE singularity  as a \emph{rational double point of type ADE}.

If we work over an algebraically closed field $k$ of characteristic $p>5$, the rational double points are realized as invariant subrings of finite subgroups $G\subset\SL_2(k)$ called \emph{binary polyhedral groups} acting linearly on the power series ring $k\llbracket u,v\rrbracket$ (\cite[\S6]{LW}, \cite{Has15}). Our aim in this section is to present such a description in the mixed characteristic case.

Let us first define binary polyhedral groups in our setting, using generators explicitly.

\begin{lemma}
\label{lem:zetaWk}
If $n\geq 1$ is an integer such that $p\nmid n$, then $W(k)$ contains a primitive $n$-th root of unity.
\end{lemma}

\begin{proof}
The algebraically closed $k$ contains a primitive $n$-th root of unity $\zeta_n$. Since $p\nmid n$, there exists a unique lift $\wt{\zeta_n}\in W(k)$ of $\zeta_n$ by Hensel's lemma. Then we can see that $\wt{\zeta_n}$ is a primitive $n$-th root of unity. (In fact, the lift coincides with the Teichm\"{u}ller lift $[\zeta_n]$ of $\zeta_n$.)
\end{proof}

\begin{definition}
We define the following subgroup schemes of $\SL_2$.
  \begin{enumerate}
  \item For an integer $n\geq 1$, we regard $\mu_{n+1}$ as the subgroup scheme of $\SL_2$ given by
  \[
  \mu_{n+1} =
  \left\{
  \begin{pmatrix}
  a & 0 \\
  0 & a^{-1}
  \end{pmatrix}
  \relmiddle| a\in \mu_{n+1}\right\} \subset \SL_2.
  \]
  \item For an integer $m\geq 2$, we define the \emph{binary dihedral group $D_{4m}$ of order $4m$} by
  \[
  D_{4m} \coloneqq \left\langle \mu_{2m}, \begin{pmatrix} 0 & -1 \\ 1 & 0 \end{pmatrix} \right\rangle \subset \SL_2.
  \]
  \item Assume that $p\geq 3$. Then $W(k)$ contains a primitive $8$-th root of unity $\zeta_8$ (\cref{lem:zetaWk}), and $\sqrt{2}\coloneqq \zeta_8+\zeta_8^{-1}\in W(k)$ satisfies $(\sqrt{2})^2=2\in W(k)^\times$. We define the \emph{binary tetrahedral group} $T$ by
  \[
  T \coloneqq \left\langle D_8, \frac{1}{\sqrt{2}} \begin{pmatrix} \zeta_8^7 & \zeta_8^7 \\ \zeta_8^5 & \zeta_8 \end{pmatrix} \right\rangle \subset \SL_2.
  \]
  \item Assume that $p\geq 3$. Then we define the \emph{binary octahedral group} $O$ by
  \[
  O\coloneqq\langle T,\mu_8\rangle\subset\SL_2.
  \]
  \item Assume that $p\neq 5$. Then $W(k)$ contains a primitive $5$-th root of unity $\zeta_5$ (\cref{lem:zetaWk}), and $\zeta_5^2-\zeta_5^3=\zeta_5^2(1-\zeta_5)\in W(k)^\times$. We define the \emph{binary icosahedral group} $I$ by
  \[
  I \coloneqq \left\langle \mu_{10},\begin{pmatrix} 0 & 1 \\ -1 & 0 \end{pmatrix}, \frac{1}{\zeta_5^2-\zeta_5^3} \begin{pmatrix} \zeta_5+\zeta_5^{-1} & 1 \\ 1 & -(\zeta_5+\zeta_5^{-1}) \end{pmatrix} \right\rangle \subset\SL_2.
  \]
  \end{enumerate}
\end{definition}

We define binary polyhedral groups $G$ under the assumption on $p$ so that $G$ is linearly reductive:

\begin{proposition-definition}
\label{prop-def:binary}
Each of the following group schemes $G$ is finite and linearly reductive over $W(k)$ and admits a functorial Reynolds operator under the indicated assumptions on $p$:
\begin{enumerate}
\item $G=\mu_{n+1}$ for some integer $n\geq1$, with no restriction on $p$.
\item $G=D_{4m}$ for some integer $m\geq2$, with $p\geq3$.
\item $G=T$ or $O$, with $p\geq5$.
\item $G=I$, with $p\geq7$.
\end{enumerate}
We refer to these group schemes as the \emph{binary polyhedral group schemes (over $W(k)$)}.
\end{proposition-definition}

\begin{proof}
(1) follows from \cref{ex:linred} (2).

(2) Let $\tau\coloneqq\begin{psmatrix}0 & -1 \\ 1 & 0 \end{psmatrix}$. Then $D_{4m}=\mu_{2m}\coprod\tau\mu_{2m}$, and thus it is finite flat. Moreover, since $p\geq 3$, $G=D_{4m}$ admits a functorial Reynolds operator by the exact sequence $1\to \mu_{2m} \to D_{4m} \to \Z/2\Z\to 0$ (\cref{ex:linred} (3)).

(3) Since $p\geq 3$, both $T$ and $O$ are the constant group schemes associated with finite subgroups $\SL_2(W(k))$ of order $24$ and $48$, respectively. Hence, if $p\geq 5$, then both $T$ and $O$ admit funtorial Reynolds operators (\cref{ex:linred} (1)).

(4) can be shown similarly to (3).
\end{proof}

We are now ready to state our first main theorem.

\begin{theorem}
\label{thm:RDPG}
Let $R=S/(f)$ be a rational double point of type ADE, where $(S,\n,k)$ is a complete regular local ring of mixed characteristic $(0,p)$ and $f\in\n$. 
Then there exists a choice, depending on the ADE type of $R$, of an isomorphism
\[
S\cong W(k)\llbracket x,y,z\rrbracket/(p-Q(x,y,z))
\]
for some $Q(x,y,z)\in (x,y,z)$, a binary polyhedral group $G\subset\SL_2$, and $G$-invariant homogeneous polynomials $X(u,v),Y(u,v),Z(u,v)\in W(k)[u,v]$ of degree at least two such that $R\cong A^G$, where 
\[
A\coloneqq W(k)\llbracket u,v\rrbracket/(p-Q(X(u,v),Y(u,v),Z(u,v)))
\]
is a ramified two-dimensional complete regular local ring endowed with the $G$-action induced by its linear action on $W(k)\llbracket u,v\rrbracket$. 
More precisely, for each ADE type, the group scheme $G$ may be chosen as in \textup{Table~\ref{tb:RDP}}.
\end{theorem}

\begin{proof}
We divide the proof into cases according to the type of $R$. In each case, we make a change of variables and replace the equation $f$ as in Table~\ref{tb:RDP} by a more convenient equation. The same change of variables also transforms $Q$, and we continue to denote the resulting power series by $Q$.

\textbf{(The case of type $A_n$)}: Since $p\neq 2$, $W(k)$ contains a primitive $4$-th root of unity $i$. Replacing $x$ and $y$ by $x+iy$ and $x-iy$, respectively, we may assume that $f=xy+z^{n+1}$. We denote again by $Q$ the power series applied the same change of variables to $Q$ .
Now a $\mu_{n+1}$-algebra is nothing but a $\Z/(n+1)\Z$-graded $W(k)$-algebra. The formal power series ring $W(k)\llbracket u,v\rrbracket$ is a $\mu_{n+1}$-algebra with $\deg(u)=1$ and $\deg(v)=-1$. Then it is easy to see that $x\mapsto u^{n+1}$, $y\mapsto -v^{n+1}$, and $z\mapsto uv$ give the isomorphism of $W(k)$-algebras
\[
W(k)\llbracket x,y,z\rrbracket/(xy+z^{n+1}) \xr{\cong} W(k)\llbracket u^{n+1},v^{n+1},uv\rrbracket = W(k)\llbracket u,v\rrbracket^{\mu_{n+1}}.
\]
Taking modulo $p-Q(x,y,z)$, we obtain the isomorphism
\[
R\cong W(k)\llbracket x,y,z\rrbracket/(p-Q(x,y,z),xy+z^{n+1}) \xr{\cong} W(k)\llbracket u,v\rrbracket^{\mu_{n+1}}/(p-Q(u^{n+1},-v^{n+1},uv))=A^{\mu_{n+1}},
\]
where $A\coloneqq W(k)\llbracket u,v\rrbracket/(p-Q(u^{n+1},-v^{n+1},uv))$.

\textbf{(The case of type $D_n$)}:
As in the case of type $A_n$, one can mimic the proof of \cite[Corollary 3.10]{Has15}. Alternatively, one can reduce to the case over $k$ as follows: since $p\geq 3$, we may assume that $f=-x^2+y^2z+4(-z)^{n-1}$. Set
\begin{equation*}
\left\{ \,
    \begin{aligned}
    X(u,v) &\coloneqq uv(u^{2(n-2)}-(-1)^{n-2}v^{2(n-2)}), \\
    Y(u,v) &\coloneqq u^{2(n-2)}+(-1)^nv^{2(n-2)}, \\
    Z(u,v) &\coloneqq u^2v^2,
    \end{aligned}
\right.
\end{equation*}
and $A\coloneqq W(k)\llbracket u,v\rrbracket/(p-Q(X,Y,Z))$.
Since $X,Y,Z\in A^{D_{4(n-2)}}$ and $Y^2Z+4(-Z)^{n-1}=X^2$, we have a homomorphism of $W(k)$-algebras
\[
R\cong W(k)\llbracket x,y,z\rrbracket/(p-Q(x,y,z),f) \to A^{D_{4(n-2)}}
\]
sending $x,y,z$ to $X,Y,Z$, respectively. Since both $R$ and $A^{D_{4(n-2)}}$ are two-dimensional domains, we only have to prove the surjectivity. Since both $R$ and $A^{D_{4(n-2)}}$ are $p$-adically complete, it suffices to check the surjectivity after taking modulo $p$ (\cite[Theorem 8.4]{Mat2}). Due to \cref{lem:AGpure}, the map obtained by taking modulo $p$ is
\[
k\llbracket x,y,z\rrbracket/(\ol{Q}(x,y,z),\ol{f}) \to \left( k\llbracket u,v\rrbracket/(\ol{Q}(X,Y,Z))\right)^{D_{4(n-2)}}.
\]
By \cref{lem:AGpure} again, this map is the same as
\[
k\llbracket x,y,z\rrbracket/(\ol{Q}(x,y,z),\ol{f}) \to k\llbracket u,v\rrbracket^{D_{4(n-2)}}/(\ol{Q}(X,Y,Z)).
\]
Hence we are reduced to show the case over $k$, that is, to establish the isomorphism
\[
k\llbracket x,y,z\rrbracket/(\ol{f}) \xr{\cong} k\llbracket u,v\rrbracket^{D_{4(n-2)}},
\]
which follows from the proof of \cite[Corollary 3.10]{Has15}.

\textbf{(The case of type $E_6$)}: Since $p\geq 5$, we may assume that $f=-x^2+y^3+108z^4$. Set
\begin{equation*}
\left\{ \,
    \begin{aligned}
    X_T(u,v) &\coloneqq u^{12}-33u^8v^4-33u^4v^8+v^{12},\\
    Y_T(u,v) &\coloneqq u^8+14u^4v^4+v^8,\\
    Z_T(u,v) &\coloneqq uv(u^4-v^4),
    \end{aligned}
\right.
\end{equation*}
and $A\coloneqq W(k)\llbracket u,v\rrbracket/(p-Q(X_T,Y_T,Z_T))$.
Then, by an argument similar to that in the case of type $D_n$, one can show the isomorphism
\[
R\cong W(k)\llbracket x,y,z\rrbracket/(p-Q(x,y,z),f) \xr{\cong} A^T.
\]
See \cite[6.15]{LW} in the case over $k$.\footnote{We denote the polynomials $X_G,Y_G,Z_G$ and the variables $x,y,z$ in \cite{LW} by $Y_G,Z_G,X_G$ and $y,z,x$, respectively.}

\textbf{(The case of type $E_7$)}: Since $p\geq 5$, we may assume that $f=x^2+y(108y^2-z^3)$. Set
\begin{equation*}
\left\{ \,
    \begin{aligned}
    X_O(u,v) &\coloneqq X_T(u,v)Z_T(u,v) = uv(u^4-v^4)(u^{12}-33u^8v^4-33u^4v^8+v^{12}), \\
    Y_O(u,v) &\coloneqq Z_T(u,v)^2 = (u^5v-uv^5)^2, \\
    Z_O(u,v) &\coloneqq Y_T(u,v) = u^8+14u^4v^4+v^8, \\
    \end{aligned}
\right.
\end{equation*}
and $A\coloneqq W(k)\llbracket u,v\rrbracket/(p-Q(X_O,Y_O,Z_O))$.
Then, by an argument similar to that in the case of type $D_n$, one can show the isomorphism
\[
R\cong W(k)\llbracket x,y,z\rrbracket/(p-Q(x,y,z),f) \xr{\cong} A^O.
\]
See \cite[6.16]{LW} in the case over $k$.

\textbf{(The case of type $E_8$)}: Since $p\geq 7$, we may assume that $f=-x^2+y^3+1728z^5$. Set
\begin{equation*}
\left\{ \,
    \begin{aligned}
    X_I(u,v) &\coloneqq (u^{30}+v^{30})+522(u^{25}v^5-u^5v^{25})-10005(u^{20}v^{10}+u^{10}v^{20}), \\
    Y_I(u,v) &\coloneqq (u^{20}+v^{20})-228(u^{15}v^5-u^5v^{15})+494u^{10}v^{10}, \\
    Z_I(u,v) &\coloneqq uv(u^{10}+11u^5v^5+v^{10}), \\
    \end{aligned}
\right.
\end{equation*}
and $A\coloneqq W(k)\llbracket u,v\rrbracket/(p-Q(X_I,Y_I,Z_I))$.
Then, by an argument similar to that in the case of type $D_n$, one can show the isomorphism
\[
R\cong W(k)\llbracket x,y,z\rrbracket/(p-Q(x,y,z),f) \xr{\cong} A^I.
\]
See \cite[6.17]{LW} in the case over $k$.
\end{proof}

\section{Frobenius lifts on invariant subrings}
\label{s:FrobG}

Throughout this section, we fix a prime number $p$; $R$ is a ring equipped with a Frobenius lift $\varphi_R\colon R\to R$ and $G$ is a linearly reductive group scheme over $R$. Here, by a \emph{Frobenius lift}, we mean a ring endomorphism $\varphi_R$ such that the induced endomorphism of $R/pR$ coincides with the absolute Frobenius. 
We are going to investigate conditions under which the invariant subring $A^G$ admits a Frobenius lift, assuming that the $G$-algebra $A$ admits one.

Now let us state an immediate but important lemma.

\begin{lemma}
\label{lem:Grest}
Let $A$ be a $G$-algebra. If a Frobenius lift $\varphi_A$ satisfies $\varphi_A(A^G)\subset A^G$, then it restricts to a Frobenius lift on $A^G$.
\end{lemma}

\begin{proof}
The homomorphism $A^G/pA^G\to A/pA$ induced by the inclusion $A^G\hookrightarrow A$ is injective by \cref{lem:AGpure}. Therefore, since $\varphi_A$ is a Frobenius lift, so is its restriction $A^G\to A^G$.
\end{proof}

Since we work over the base ring $R$, it is more natural to consider lifts of the relative Frobenius, as defined below.

\begin{definition}
\label{def:relFlift}
Let $S\coloneqq\Spec(R)$, $F_S\coloneqq\Spec(\varphi_R)\colon S\to S$, and $X$ an $S$-scheme.
  \begin{enumerate}
  \item We define the $S$-scheme $X^{(p)}$ by the cartesian square
  \[
  \begin{tikzcd}
  X^{(p)} \rar \dar \ar[rd,phantom,"\scriptstyle\Box"] & X \dar \\
  S \rar["F_S"'] & S.
  \end{tikzcd}
  \]
  If $X=\Spec(A)$ is affine, we write $A^{(p)}$ for the coordinate ring of $X^{(p)}$, so that $X^{(p)}=\Spec(A^{(p)})$.
  \item A \emph{relative Frobenius lift of $X/S$} is a morphism of $S$-schemes $F_{X/S}\colon X\to X^{(p)}$ such that the induced morphism of $S_0$-schemes $X_0\to X_0^{(p)}$ coincides with the relative Frobenius of $X_0/S_0$. Here $(-)_0\coloneqq-\otimes_R(R/pR)$.
  If $X=\Spec(A)$ is affine, we write $\varphi_{A/R}\colon A^{(p)}\to A$ for the morphism of $R$-algebras corresponding to $F_{X/S}$ and call it the \emph{relative Frobenius lift of $A/R$}.
  \end{enumerate}
\end{definition}

Note that giving a relative Frobenius lift of $X/S$ amounts to the same as giving a Frobenius lift $F_X\colon X\to X$ on the scheme $X$ such that the diagram of schemes
\[
\begin{tikzcd}
X \rar["F_X"] \dar & X \dar \\
S \rar["F_S"] & S
\end{tikzcd}
\]
commutes. One of the most typical examples is the endomorphism of $R$-algebras $\varphi_{A/R}$ on the polynomial ring $A=R[x_1,\ldots,x_n]$ such that $\varphi_{A/R}(x_i)=x_i^p$ ($1\leq{}^\forall i\leq n$). As for group schemes, we have the following examples.

\begin{example}
\label{ex:FG}
  \begin{enumerate}
  \item Let $G=\ul{\Gamma}_R$ be the constant group scheme associated with a finite group $\Gamma$. Then $\left(\prod\limits_\Gamma R\right)\otimes_{R,\varphi_R}R\xr{\cong} \prod\limits_\Gamma R$, and thus one has the canonical isomorphism $G\xr{\cong} G^{(p)}$, which gives a relative Frobenius lift.
  \item Let $G=D_R(\Delta)$ be the diagonalizable group scheme associated with an abelian group $\Delta$. Then $R\Delta\otimes_{R,\varphi_R}R \xr{\cong} R\Delta$, and thus one has the canonical isomorphism $G\xr{\cong}G^{(p)}$. Hence we have the morphism $F_{G/R}\colon G\to G^{(p)}$ induced by $\Delta\to\Delta$; $d\mapsto pd$. We see that $F_{G/R}$ is a relative Frobenius lift.
  \end{enumerate}
\end{example}

We present a useful sufficient condition for a relative Frobenius lift of $A/R$ to induce a Frobenius lift on $A^G$.

\begin{proposition}
\label{prop:Fequiv}
Let $A$ be a $G$-algebra equipped with a relative Frobenius lift $\varphi_{A/R}$.
Suppose that there exists a morphism of $R$-group schemes $F\colon G\to G^{(p)}$ such that $\varphi_{A/R}$ is \emph{$F$-equivalent} in the sense that the diagram
\begin{equation}
\label{eq:psi}
\begin{tikzcd}
A^{(p)} \rar["\omega_A^{(p)}"] \dar["\varphi_{A/R}"'] & A^{(p)}\otimes_R R[G]^{(p)} \dar["\varphi_{A/R}\otimes_R F^\#"] \\
A \rar["\omega_A"'] & A\otimes_RR[G]
\end{tikzcd}
\end{equation}
commutes, where $F^\#\colon R[G]^{(p)}\to R[G]$ is the morphism of $R$-algebras corresponding to $F$. Then the Frobenius lift $\varphi_A$ associated with $\varphi_{A/R}$ restricts to a Frobenius lift on $A^G$.
\end{proposition}

\begin{proof}
Pick $x\in A^G$. Then $\omega_A(x)=x\otimes 1$, and so $\omega_A^{(p)}(x^{(p)}) = x^{(p)}\otimes 1$ in $A^{(p)}\otimes_RR[G]^{(p)}$. where $x^{(p)}\coloneqq x\otimes 1 \in A^{(p)}=A\otimes_{R,\varphi}R$. By assumption, one has
\begin{align*}
\omega_A\left(\varphi_A(x)\right) &= \omega_A\left(\varphi_{A/R}(x^{(p)})\right) \\
&= (\varphi_{A/R}\otimes_RF^\#)\left(\omega_A^{(p)}(x^{(p)})\right) \\
&= (\varphi_{A/R}\otimes_RF^\#)(x^{(p)}\otimes1) \\
&= \varphi_A(x)\otimes 1.
\end{align*}
Thus $\varphi_A(x)\in A^G$, and hence $\varphi_A(A^G)\subset A^G$. The assertion now follows from \cref{lem:Grest}.
\end{proof}


\begin{example}[Diagonalizable group actions]
\label{ex:diagF}
Let $G=D_R(\Delta)$ be the diagonalizable group scheme associated with an abelian group $\Delta$. Then a $G$-algebra is nothing but a $\Delta$-graded $R$-algebra $A=\bigoplus_{d\in\Delta}A_d$. Consider the relative Frobenius lift $F_{G/R}\colon G\to G^{(p)}$ as in \cref{ex:FG} (2). Then a relative Frobenius lift $\varphi_{A/R}$ makes the diagram \eqref{eq:psi} commutes if and only if the corresponding Frobenius lift $\varphi_A$ satisfies $\varphi_A(A_d)\subset A_{pd}$ for all $d\in\Delta$. In this case, $\varphi_A$ does indeed restrict to a Frobenius lift on $A^G=A_0$.
\end{example}

In particular, for polynomial rings and formal power series rings, we obtain the following.

\begin{proposition}
\label{prop:diagF}
Let $G=D_R(\Delta)$ be the diagonalizable group scheme associated with an abelian group, and fix elements $d_1,\ldots,d_n\in\Delta$. Let $A$ be either $R[x_1,\ldots,x_n]$ or $R\llbracket x_1,\ldots,x_n\rrbracket$, and let $G$ act trivially on $R$ and diagonally on $A$ with $x_i$ of weight $d_i$.
Let $\varphi_A\colon A\to A$ be the Frobenius lift defined by $\varphi_A|_R=\varphi_R$ and $\varphi_A(x_i)=x_i^p$ ($1\leq{}^\forall i\leq n$). Then $\varphi_A$ restricts to a Frobenius lift on $A^G=A_0$.
\end{proposition}

Next we give a construction of an $F_{G/R}$-equivariant Frobenius lift.

\begin{theorem}
\label{thm:GequivF}
Suppose that $G$ admits a functorial Reynolds operator and a relative Frobenius lift $F_{G/R}\colon G\to G^{(p)}$.
Set the polynomial ring $A\coloneqq R[x_1,\ldots,x_n]$ equipped with a linear action of $G$, by which we mean that $M\coloneqq Rx_1\oplus\cdots\oplus Rx_n$ is $G$-stable.
Then there exists an $F_{G/R}$-equivariant relative Frobenius lift $\varphi_{A/R,G}\colon A\to A^{(p)}$. In particular, $A^G$ admits a Frobenius lift.

Moreover, the similar assertion holds for $A=R\llbracket x_1,\ldots,x_n\rrbracket$ if $R$ is $p$-adically complete. 
\end{theorem}

\begin{proof}
We equip the twist $M^{(p)}\coloneqq M\otimes_{R,\varphi_R}R=Rx_1^{(p)}\oplus\cdots\oplus Rx_n^{(p)}$ with the $G$-action by restriction along $F_{G/R}$, and the $G$-module $\Hom_R(M^{(p)},A)$ (\cref{rem:HomG}).
Let $f\colon M^{(p)}\to A$ be the $R$-linear map defined by $f(x_i^{(p)})\coloneqq x_i^p$.
Then $\ol{f}\coloneqq f\ \mathrm{mod}\ p \colon (M/pM)^{(p)}\to A/pA$ is the restriction of the relative Frobenius $(A/pA)^{(p)}\to A/pA$. In particular, $\ol{f}$ is $G_0$-equivariant, where $G_0\coloneqq G\otimes_R(R/pR)$. This amounts to saying that $\ol{f}\in\Hom_{R/pR}((M/pM)^{(p)},A/pA)$ is $G_0$-invariant by \cref{rem:HomG}. Hence the functoriality of Reynolds operators $\rho$ implies that
\begin{equation}
\label{eq:rhof}
\rho_{\Hom_R(M^{(p)},A)}(f)\equiv f\ \mathrm{mod}\ p
\end{equation}
Now let $\varphi_{A/R,G}\colon A^{(p)}\to A$ be the morphism of $R$-algebras such that
\[
\varphi_{A/R,G}(x_i) \coloneqq \rho_{\Hom_R(M^{(p)},A)}(f)(x_i)\quad (1\leq{}^\forall i\leq n).
\]
Then it follows from \eqref{eq:rhof} that $\varphi_{A/R,G}$ is a relative Frobenius lift. Moreover, since $\rho_{\Hom_R(M^{(p)},A)}(f)\in\Hom_R(M^{(p)},A)^G$ is $G$-equivariant by \cref{rem:HomG}, we see that $\varphi_{A/R,G}\colon A^{(p)}\to A$ commutes the diagram \eqref{eq:psi} for $F_{G/R}$.

In the case $A=R\llbracket x_1,\ldots,x_n\rrbracket$, we only have to remark that $\varphi_{A/R,G}$ is well-defined by the $(p,x_1,\ldots,x_n)$-adic completeness of $A$.
\end{proof}

\begin{example}
Let $G=\ul{\Gamma}_R$ be the constant group scheme associated with a finite subgroup $\Gamma\subset\GL_n(R)$ such that $\#\Gamma\in R^\times$.
Set $A\coloneqq R[x_1,\ldots,x_n]$, and let $G$ act linearly on $A$ from the left. Let $\varphi_A\colon A\to A$ be the Frobenius lift defined by $\varphi_A|_R=\varphi_R$ and $\varphi_A(x_i)=x_i^p$ ($1\leq{}^\forall i\leq n$). Then the ring homomorphism $\varphi_{A,\Gamma}\colon A\to A$ defined by $\varphi_{A,\Gamma}|_R = \varphi_R$ and
\[
\varphi_{A,\Gamma}(x_i)\coloneqq \frac{1}{\#\Gamma}\sum_{\gamma\in \Gamma}\gamma\cdot \varphi_A(\gamma^{-1}\cdot x_i)\quad (1\leq{}^\forall i\leq n)
\]
is a $G$-equivariant Frobenius lift on $A$.
\end{example}

\section{Applications to rational double points and perfectoid towers}
\label{s:RDPtower}

Throughout this section, $k$ is an algebraically closed field of characteristic $p>0$ and $W(k)$ is the ring of Witt vectors. All group schemes are over $W(k)$.

\begin{proposition}
\label{prop:binaryF}
For every binary polyhedral group scheme $G$ over $W(k)$ \textup{(\cref{prop-def:binary})}, there exists a relative Frobenius lift $F_{G/W(k)}\colon G\to G^{(p)}$.
\end{proposition}

\begin{proof}
By \cref{ex:FG}, we only have to show the case $G=D_{4m}$. 
Recall that $G=\mu_{2m}\coprod \tau\mu_{2m}$, where $\tau\coloneqq\begin{psmatrix}0 & -1 \\ 1 & 0 \end{psmatrix}$. We define the morphism of $W(k)$-schemes $F_{G/W(k)}\colon G\to G^{(p)}$ on the two components by
\[
F_{G/W(k)}(\sigma)=\left(\sigma^p\right)^{(p)},\quad F_{G/W(k)}(\tau\sigma)=\left(\tau\sigma^p\right)^{(p)}\quad ({}^\forall \sigma\in\mu_{2m}).
\]
To check that it is a group homomorphism, it suffices to verify the compatibility with multiplication on the four components of
\[
G\times G = \left( \mu_{2m}\times\mu_{2m}\right) \coprod \left( \mu_{2m}\times \tau\mu_{2m}\right) \coprod \left( \tau\mu_{2m}\times\mu_{2m}\right) \coprod \left(\tau\mu_{2m}\times\tau\mu_{2m}\right).
\]
Pick $\sigma_1,\sigma_2\in \mu_{2m}$. We have $F_{G/W(k)}(\sigma_1\sigma_2) = \left((\sigma_1\sigma_2)^p\right)^{(p)} = \left(\sigma_1^p\right)^{(p)} \left(\sigma_2^p\right)^{(p)} = F_{G/W(k)}(\sigma_1\sigma_2)$.
Since $\tau^{-1}\sigma_1\tau=\sigma_1^{-1}$, we also have
\[
F_{G/W(k)}(\sigma_1(\tau\sigma_2)) = F_{G/W(k)}(\tau\sigma_1^{-1}\sigma_2) = \left(\tau \sigma_1^{-p}\sigma_2^p\right)^{(p)} = \left(\sigma_1^p\tau\sigma_2^p\right)^{(p)} = F_{G/W(k)}(\sigma_1)F_{G/W(k)}(\tau\sigma_2).
\]
Similarly, we have $F_{G/W(k)}(\sigma_1(\tau\sigma_2)) = F_{G/W(k)}(\tau\sigma_1)F_{G/W(k)}(\sigma_2)$. Moreover, since $\tau^2=-1$, we have
\begin{align*}
F_{G/W(k)}((\tau\sigma_1)(\tau\sigma_2)) &= F_{G/W(k)}((-1)\sigma_1^{-1}\sigma_2) \\
&= \left((-1)^p\sigma_1^{-p}\sigma_2^{-p}\right)^{(p)} \\
&\stackrel{p\geq 3}{=} \left((-1)\sigma_1^{-p}\sigma_2^{-p}\right)^{(p)} \\
&= \left(\tau\sigma_1^p\right)^{(p)} \left(\tau\sigma_2^p\right)^{(p)} \\
&= F_{G/W(k)}(\tau\sigma_1)F_{G/W(k)}(\tau\sigma_2).
\end{align*}
Hence $F_{G/W(k)}$ is a homomorphism of group schemes. Finally, $F_{G/W(k)}\ \mathrm{mod}\ p$ is given by raising each matrix entry to its $p$-th power. Therefore, $F_{G/W(k)}$ is the relative Frobenius lift.
\end{proof}

By \cref{prop:binaryF,thm:GequivF}, we have the following result.

\begin{corollary}
\label{cor:RDPphi}
Let $G$ be a binary polyhedral group scheme over $W(k)$, and consider its linear action on $W(k)\llbracket u,v\rrbracket$. Then the invariant subring $W(k)\llbracket u,v\rrbracket^G$ admits a relative Frobenius lift.
\end{corollary}

\begin{remark}
In the case $G$ is either $\mu_{n+1}$ or $D_{4m}$, we may apply \cref{prop:Fequiv} rather than \cref{thm:GequivF}. Indeed, we claim that the relative Frobenius lift $\varphi_{/W(k)} \colon W(k)\llbracket u,v\rrbracket \to W(k)\llbracket u,v\rrbracket$ defined by $\varphi_{/W(k)}(u)=u^p$ and $\varphi_{/W(k)}(v)=v^p$ is $F_{G/R}$-equivalent, where $F_{G/R}\colon G\to G^{(p)}$ is the relative Frobenius lift as in \cref{ex:FG} (2) or the proof of \cref{prop:binaryF}.

The case $G=\mu_{n+1}$ follows from \cref{prop:diagF}.
Consider the case $G=D_{4m}=\mu_{2m}\coprod\tau\mu_{2m}$, where $\tau\coloneqq\begin{psmatrix}0 & -1 \\ 1 & 0 \end{psmatrix}$. The action of $\mu_{2m}$ on $W(k)\llbracket u,v\rrbracket$ induces a ($\Z/2m\Z$)-grading for which $u$ and $v$ have degrees $1$ and $-1$, respectively, while $\tau\cdot u=-v$ and $\tau\cdot v=u$.
For the action of $\mu_{2m}$, it follows from \cref{ex:diagF} that $\varphi_{/W(k)}$ is $F_{G/W(k)}|_{\mu_{2m}}$-equivariant.
Moreover, we have
\begin{align*}
\varphi_{/W(k)}(\tau\cdot u) &=\varphi_{/W(k)}(-v)=-v^p=\tau\cdot u^p=\tau\cdot \varphi_{/W(k)}(u), \\
\varphi_{/W(k)}(\tau\cdot v) &=\varphi_{/W(k)}(u)=u^p=\tau\cdot\varphi_{/W(k)}(v).
\end{align*}
Hence $\varphi_{/W(k)}$ is $F_{G/R}$-equivariant.
\end{remark}

Finally, we prove that the local ring of a rational double point admits a perfectoid tower. We refer to \cite[Definition 3.21]{INS25} for the definition of perfectoid towers. For our purposes, we only use Ishizuka's construction of perfectoid towers from prisms given in \cite[Theorem 4.5]{Ish26}. Here a \emph{prism} is a pair consisting of a ring endowed with a \emph{$\delta$-structure} and an ideal satisfying certain conditions; typically, a $\delta$-structure is given by a Frobenius lift. See \cite[\S\S2 and 3]{BS22} for details. 

\begin{theorem}
\label{thm:RDPtower}
Let $R=S/(f)$ be a rational double point of type ADE, where $(S,\n,k)$ is a three-dimensional complete regular local ring of mixed characteristic $(0,p)$ with algebraically closed residue field $k$ and $f\in\n$.
Choose an isomorphism $S\cong W(k)\llbracket x,y,z\rrbracket/(p-Q(x,y,z))$, a binary polyhedral group $G\subset \SL_{2,W(k)}$, and $G$-invariant polynomials $X(u,v),Y(u,v),Z(u,v)\in W(k)[u,v]$ as in \textup{\cref{thm:RDPG}}, so that
\[
R\cong A^G,\quad A\coloneqq W(k)\llbracket u,v\rrbracket/(p-Q(X(u,v),Y(u,v),Z(u,v))).
\]
Equip $W(k)\llbracket u,v\rrbracket^G$ with the $\delta$-structure induced by the Frobenius lift $\varphi$ \textup{(\cref{cor:RDPphi})}, and set $I\coloneqq (p-Q(X(u,v),Y(u,v),Z(u,v))) \subset W(k)\llbracket u,v\rrbracket^G$.
  \begin{enumerate}
  \item $\left(W(k)\llbracket u,v\rrbracket^G, I\right)$ is an orientable prism.
  \item The inductive system of rings
  \[
  R\cong \frac{W(k)\llbracket u,v\rrbracket^G}{I} \xr{\varphi} \frac{W(k)\llbracket u,v\rrbracket^G}{\varphi(I)W(k)\llbracket u,v\rrbracket^G} \xr{\varphi}\cdots \xr{\varphi} \frac{W(k)\llbracket u,v\rrbracket^G}{\varphi^i(I)W(k)\llbracket u,v\rrbracket^G} \xr{\varphi} \cdots
  \]
  is a perfectoid tower arising from the pair $(R,(p))$. Its tilt $($\cite[Definition 3.29]{INS25}$)$ is isomorphic to the inductive system of $\F_p$-algebras
  \[
  k\llbracket u,v\rrbracket^{G_k} \to k\llbracket u,v\rrbracket^{G_k} \to \cdots \to k\llbracket u,v\rrbracket^{G_k} \to \cdots
  \]
  where $G_k\coloneqq G\otimes_{W(k)}k$ and transition maps are given by the absolute Frobenius.
  \end{enumerate}
\end{theorem}

\begin{proof}
(1) Since $W(k)\llbracket u,v\rrbracket^G$ is a complete (Noetherian) local ring with $p,I$ contained in the maximal ideal $\m$, it suffices to show that
\[
d\coloneqq p-Q(X(u,v),Y(u,v),Z(u,v))
\]
is a distinguished element (\cite[Lemma 2.25]{BS22}). Let $\delta\colon W(k)\llbracket u,v\rrbracket^G\to W(k)\llbracket u,v\rrbracket^G$ denote the $\delta$-map induced by the Frobenius lift $\varphi$ on $W(k)\llbracket u,v\rrbracket^G$. Since $X(u,v)$, $Y(u,v)$, and $Z(u,v)$ are homogeneous of degree at least two, one has $Q(X(u,v),Y(u,v),Z(u,v))\in \m^2$, and thus $\delta(Q(X(u,v),Y(u,v),Z(u,v))) \in \m$ by \cite[Lemma 3.3]{IN26}. Moreover, $\delta(p)$ is a unit. Then a direct calculation of $\delta(d)$ shows that $d$ is a distinguished element, as desired.

(2) Since $W(k)\llbracket u,v\rrbracket^G/(p) \cong k\llbracket u,v\rrbracket^{G_k}$ by \cref{lem:AGpure}, we only have to apply \cite[Theorem~4.5]{Ish26}. For this we need to show the following.

\begin{claim}
  \begin{enumroman}
  \item $p,d$ is a regular sequence on $W(k)\llbracket u,v\rrbracket^G$.
  \item $W(k)\llbracket u,v\rrbracket^G/(p)$ is $p$-root closed in $(W(k)\llbracket u,v\rrbracket^G/(p))[1/d]$, that is, if $x\in (W(k)\llbracket u,v\rrbracket^G/(p))[1/d]$ satisfies $x^p\in W(k)\llbracket u,v\rrbracket^G/(p)$, then $x\in W(k)\llbracket u,v\rrbracket^G/(p)$.
  \end{enumroman}
\end{claim}

\begin{pfclaim}
(i) follows from the fact that both $W(k)\llbracket u,v\rrbracket^G$ and $W(k)\llbracket u,v\rrbracket^G/(p)\cong k\llbracket u,v\rrbracket^{G_k}$ are integral domains.

(ii) follows from the fact that $W(k)\llbracket u,v\rrbracket^G/(p)\cong k\llbracket u,v\rrbracket^{G_k}$ is a normal domain by \cite[Proposition 6.4.1]{BH}.
\end{pfclaim}

This completes the proof.
\end{proof}

Note that \cref{thm:RDPtower} (2) shows that the small tilt $R^{s.\flat}$ of $R$ is isomorphic to $k\llbracket u,v\rrbracket^{G_k}$. That is, under the small tilting operation $R\leadsto R^{s.\flat}$, rational double points of ADE type over $W(k)$ correspond to these over $k$, with the ADE type preserved. This observation further supports the potential applicability of the theory of perfectoid towers to singularity theory in mixed characteristic.


\end{document}